\documentstyle{amsart}

\newcommand{\nc}{\newcommand}
\renewcommand{\u}{{\hbox{\bf u}}}
\renewcommand{\b}{{\hbox{\bf b}}}
\nc{\C}{{\cal C}}
\renewcommand{\O}{{\cal O}}
\nc{\St}{\operatorname{St}^{\bullet}}
\nc{\B}{{\cal B}}
\nc{\N}{{\cal N}}
\nc{\gr}{\operatorname{gr}}
\nc{\g}{\frak g}
\nc{\Ind}{\operatorname{Ind}}
\nc{\Coind}{\operatorname{Coind}}
\nc{\Ker}{\operatorname{Ker}}
\nc{\Coker}{\operatorname{Coker}}
\nc{\dirlim}{\underset{\rightarrow}{\operatorname{lim}}}
\nc{\invlim}{\underset{\leftarrow}{\operatorname{lim}}}
\nc{\Ext}{\operatorname{Ext}^{\bullet}}
\nc{\ext}{\operatorname{Ext}}
\nc{\Tor}{\operatorname{Tor}_{\bullet}}
\nc{\tor}{\operatorname{Tor}}
\nc{\Tors}{\operatorname{Tor}_{\frac \infty 2+\bullet}}
\nc{\Exts}{\operatorname{Ext}^{\frac \infty 2+\bullet}}
\nc{\Hom}{\operatorname{Hom}}
\renewcommand{\mod}{\operatorname{-mod}}
\nc{\Mod}{\operatorname{Mod}}
\renewcommand{\Bar}{\operatorname{Bar}}
\nc{\Barb}{\operatorname{Bar}^{\bullet}}
\nc{\upC}{{\cal C}^{\uparrow}}
\nc{\dC}{{\cal C}^{\downarrow}}
\nc{\map}{\longrightarrow}
\nc{\Z}{{\mbox{\bf Z}}}
\nc{\Q}{{\mbox{\bf Q}}}
\nc{\bs}{\bigskip\\}
\nc{\ms}{\medskip\\}
\nc{\tilBar}{\widetilde{\operatorname{Bar}}}
\nc{\tilBarb}{\widetilde{\operatorname{Bar}}^{\bullet}}
\nc{\linBar}{\overline{\operatorname{Bar}}}
\nc{\linBarb}{\overline{\operatorname{Bar}}^{\bullet}}
\nc{\til}{\widetilde}
\nc{\oppA}{A^{\sharp}}
\nc{\Lemma}{{\bf Lemma:\ }}
\nc{\Theorem}{{\bf Theorem:}\ }
\nc{\Cor}{{\bf Corollary:}\ }
\nc{\Def}{{\bf Definition:}\ }
\nc{\Rem}{{\bf Remark:}\ }
\nc{\Prop}{{\bf Proposition:}\ }
\nc{\Con}{{\bf Conjecture:}\ }
\nc{\dok}{{\bf Proof.}\ }
\nc{\bul}{^{\bullet}}
\nc{\stand}{C^{\frac \infty 2+\bullet }}
\nc{\Cone}{\operatorname{Cone}}
\nc{\supp}{\operatorname{supp}}
\nc{\hgt}{\operatorname{ht}}
\nc{\sqbinom}{\fracwithdelims[][0pt]}
\author{S. M. Arkhipov}
\address{Independent University of Moscow, Pervomajskaya st. 16-18,
Moscow 105037, Russia}
\email{hippie@@ium.ips.ras.ru}
\thanks{Partially supported by Soros foundation.}
\title{Semiinfinite cohomology of quantum groups}

\begin{document}

\maketitle
\section{Introduction}

Semiinfinite cohomology of Lie algebras appeared in mathematics more  than 10
years ago (see [F]),  and yet it belongs to the area of homological  algebra
existing partly in the form of folklore. A remarkable breakthrough was achieved
by A.~Voronov (see  [V]) who managed to define the semiinfinite cohomology in
the derived category setting.  Yet the general definition of semiinfinite
cohomology of associative algebras has been unknown.

\subsection{}
One of the aims of this paper is to give a
rigorous construction
of the functor of semiinfinite cohomology of arbitrary (graded) associative
algebras
defined in the corresponding derived categories of graded modules.
The basic setup here includes a graded associative algebra
$A$ with two  graded subalgebras $B$
and $N$ such that $A=B\otimes N$ as a graded vector space.
These conditions are satisfied in particular in the case of the
universal enveloping algebra of
a graded Lie algebra, but the general case is much wider.
We will show that the semiinfinite cohomology of the universal
enveloping algebra
coincides with the corresponding Lie algebra semiinfinite cohomology
(see  [V], [F]).

Our construction of the standard complex for the calculation of the
semiinfinite cohomology suggests  that semiinfinite cohomology could be realized
by some combination of the standard $\tor$ and $\ext$ functors. Our
considerations do not give such a realization, yet it exists in suitably chosen
triangulated categories. This will be explained elsewhere.

\subsection{}
Let us describe the structure of the paper.
In the second section we recall several basic facts about
quantum groups at roots of unity.  In the third section we give the
definition of the semiinfinite cohomology of associative algebras and
prove some basic results about them.  In the fourth and fifth sections we
consider an example of the  algebra $A$ equal to the finite dimensional
Hopf algebra $\u$ (finite quantum group) introduced  by
G.Lusztig in [L1]. In this case the
semiinfinite homology $\Tors^{\C}(k,Y)$
appeared in [FiS]. They calculated the cohomology of the space of
configurations of points on the projective line $\hbox{\bf P}^1$  with
coefficients in the sheaf ${\cal Y}$ equal to the localization of the
module $Y$ to the point $0 \in \hbox{\bf P}^1$.

It is  known that the semisimple Lie algebra $\frak g$ acts on the
cohomology
of the trivial module over the finite quantum group with the same
Cartan matrix (see  [GK]).  In the fourth section we show that this fact
holds for  semiinfinite cohomology of the trivial module and
calculate the character of this $\g$-module.
Unfortunately even in the simplest case the representation itself remains
unknown.  B.Feigin has proposed a conjecture
describing the $\g$-module of  semiinfinite cohomology
of the trivial module over the finite quantum group in terms of
distributions on the nilpotent cone of $\g$. In Appendix~A we prove
some facts confirming the conjecture on the level of  characters.

The main result of the fifth section is   Theorem~\ref{embed}, stating
that  conformal blocks are naturally embedded into the semiinfinite
$\tor$ spaces (see  the exact statement in section 5).  This Theorem
along with the results of the papers [Fi], [FiS] sheds  light on the
conjecture of Feigin, Schechtman and Varchenko about the integral
representation of conformal blocks (see  [FSV]). Namely, the above
results imply the following statement:  the local system of conformal
blocks on the space of  configurations of points on $\hbox{\bf P}^1$
is  a direct summand in the direct image of some perverse sheaf on a
larger configuration space. The perverse sheaf itself is the
Goresky-MacPherson extension of a one-dimensional local system.  The
example~\ref{example} shows that  the local system of conformal blocks
is in general a {\it proper} direct summand in the direct image of the
above Goresky-MacPherson sheaf.

In Appendix~B we present several Voronov's results on semiinfinite
homological algebra (see [V]) and make an uttempt to realize
semiinfinite cohomology of associative algebras as a two-sided derived
functor (in spirit of [V], 3.9).

This paper grew out of attempts to understand the natural general setting
for semiinfinite cohomology. I was introduced to the subject by B.Feigin
back in 1993. He also formulated the conjectural answer for the
semiinfinite cohomology of finite quantum groups. Thus the present paper
owes its very existence to B.Feigin. I am also greatly indebted to
M.Finkelberg and L.Positselsky for many helpful discussions. I would like to
thank D.Timashev for bringing the paper [H] to my attention.

\section*{Notation}
Throughout the paper we use the following notation.

\advance\leftmargini by 50pt
\begin{itemize}
\item[$(a)_{i,j=1}^r\ $] \ \ \ \ is a Cartan matrix of the finite type
\item[$d_1,\ldots ,d_r$] \ \ \ \ $\in \{1,2,3\}$ such that
$(d_ia_{ij})$ is symmetric
\item[$R$] \ \ \ \ is the root system
corresponding to $(a_{ij})$
\item[$R^{\pm}$] \ \ \ \ is the set of
positive (resp. negative) roots \item[$\rho$] \ \ \ \ $=
{\frac 1 2}\underset{\alpha\in R^+}{\sum}\alpha$ \item[$\Sigma$] \ \ \ \
$=\{\alpha_1,\ldots ,\alpha_r\}$ is the set of simple roots
    \item[$\operatorname{ht}\beta$] \ \ \ \ $=\sum_{i=1}^r b_i,$,
    where $\beta=\sum_{i=1}^rb_i\alpha_i\in R^+$ is the height
function on the \\ \mbox{}\ \ \ \ set of the positive roots \item[$X$]
\ \ \ \ is the weight lattice of $R$ \item[$Y$] \ \ \ \ is the root
    lattice of $R$, i.~e.
$Y=\Z\alpha_1\oplus\ldots\oplus\Z\alpha_r\hookrightarrow X$
    \item[$Y^{\pm}$] \ \ \ \ is the subsemigroup in $Y$ generated by
the set $\Sigma$ (resp. by $-\Sigma$)
\item[$(\cdot |\cdot )$] \ \ \ \
    is a scalar product defined on $\Sigma\subset X$ by the formula
    $(\alpha_i|\alpha_j)=d_ia_{ij}$ \\ \mbox{}\ \ \ \ and extended to
$X$ by bilinearity
\item[$\alpha^{\vee}_1,\ldots ,\alpha^{\vee}_r$] \
\ \ \ are the simple coroots
\item[$W$] \ \ \ \ is the Weyl group
corresponding to $R$
\item[$\g$] \ \ \ \ is the semisimple Lie algebra
with the Cartan matrix $(a)_{ij}$
\item[$G$] \ \ \ \ is the simply
connected Lie group with  Lie algebra $\g$
\item[$\Q (v)$] \ \ \ \ is
the field of rational functions in the indeterminate $v$
   \item[{$[m]!_d$}] \ \ \ \ $=\prod_{j=1}^m {\frac
{v^{dj}-v^{-dj}}{v^d-v^{-d}}}\in\Q(v)$ where $m,d\in{\mbox{\bf N}}$
   \item[{$\sqbinom mt_d$}] \ \ \ \ $=\prod_{j=1}^t{\frac
{v^{d(m-j+1)}-v^{-d(m-j+1}} {v^{dj}-v^{-dj}}}\in \Q (v)$ where
$m\in\Z,\ t,d\in{\mbox{\bf N}}$ \end{itemize}

\section{Quantum groups at roots of unity}

In this section we collect several  well-known facts about finite
 quantum groups that we will need later.

\subsection{}
For every symmetrizable Cartan matrix $(a)_{i,j=1}^r$
Drinfeld and Jimbo constructed a Hopf algebra $U_1$ over the field
$\Q (v)$ of rational functions with the generators $E_i$,
$F_i$, $K_i$, $K_i^{-1}$, $i=1,\ldots ,r$, and the following relations:
$$
K_iK_j=K_jK_i,\ K_iK_i^{-1}=K_i^{-1}K_i=1
$$
$$
K_iE_j=v^{d_ia_{ij}}E_jK_i,\ K_iF_j=v^{-d_ia_{ij}}F_jK_i
$$
$$
E_iF_j-F_jE_i=\delta_{ij}{\frac {K_i-K_i^{-1}}{v^{d_i}-v{-d_i}}}
$$
$$
\sum_{r+s=1-a_{ij}}(-1)^s\sqbinom{1-a_{ij}}s_{d_i} E_i^rE_jE_i^s=0
\text{ if } i\ne j
$$
(see also  [L1]).
This algebra is
called the quantum universal enveloping algebra of the corresponding
Kac-Moody Lie algebra or  the quantum group.

\subsection{}
Given an integer  $p>1$ prime to the nonzero elements of the
fixed Cartan matrix, we  choose a primitive $p$-th root of unity
$\zeta $ and set $k=\Q (\zeta )$.  De Concini and Kac introduced
the $k$-algebra $U_2$ with generators
$E_i$, $F_i$, $F_i$, $ K_i$, $K_i^{-1}$ and relations similar to the
ones for $U_1$ but with
$v$ replaced by $\zeta$ (see [DCK], 1.5).

Denote the subalgebra   in $U_2$ generated by all $E_i$
(resp. by all $F_i$, resp. by all $K_i$ and $K_i^{-1}$) by $U_2^{+}$
(resp. by $U_2^{-}$, resp. by $U_2^0$). Then $U_2^0$ is the algebra of
Laurent polynomials in the commuting variables
$K_i$. We define the following $X$-grading on
$U_2$: $\deg E_i=\alpha _i$, $\deg F_i=-\alpha_i$, $\deg K_i=0$.

Like in  classical universal enveloping algebras,
for every $\beta \in R^{+}$ one can define the root elements
$E_\beta \in (U_2)_\beta $ and $F_{-\beta}
\in (U_2)_{-\beta }$, in particular $E_{\alpha_i}=E_i,\ F_{-\alpha_i}=F_i$.
(see   [L1]).

\subsubsection{}
\Lemma (see [DCK], 1.7) The elements
$$
\underset{\beta \in R^{+}}{\prod}E_\beta ^{n(\beta )}
\prod_{i=1}^r K_i^{s(i)}
\underset{\beta \in R^{+}}{\prod } F_{-\beta }^{m(\beta )}
$$
are linearly independent and generate
$U_2$ as a vector space.
Here $n(\beta )$ and $m(\beta )$ run through nonnegative integers
 and $s(i)$ run through arbitrary integers.$\qed $

\subsubsection{} \Cor
The multiplication in $U_2$ defines an  isomorphism of the vector
spaces $ U_2^{+}\otimes U_2^0\otimes U_2^{-}\map U_2.\qed $ \bs All
$E_\beta ^p$ and $F_{-\beta}^p,\ \beta\in R^+,$ and $K_i^p,\
i=1,\ldots ,r$, belong to the center of the algebra $U_{2}$ (see
[DCK], 3.1).  Denote the ideals generated by the elements
$\{E_{\beta}^p,F_{-\beta}^p|\beta\in R^+\}$ in the algebras
$U^{+},U^{-}$, and $U$ respectively by $I^{+}$, $I^{-}$ and $I$.
Obviously $I^{\pm }=I\cap U_2^{\pm }$.

Set $U=U_2/I$. Let $U^{\pm }$ and $U^0$ be the images of the
corresponding  subalgebras under the projection.  Note that  $U^+$
(resp. $U^-$) is  $Y^+$-graded (resp. $Y^-$-graded).  We obtain the
base in $U$ that consists of monomials similar to those of the
previous Lemma but with $0\leq m(\beta ),n(\beta )\leq p$ for every
$\beta \in R^{+}$ and arbitrary $s(i)\in\Z$.

\subsection{}
There exists an alternative version of the quantum
group  defined by Lusztig (see  [L1], 8.1).
Namely Lusztig considers the $\Q [v,v^{-1}]$-subalgebra
$U_{\Z}\subset U_1$ generated by the elements
$$
E_i^{(n)}:=E_i^n/[n]!_{d_i},\ F_i^{(n)}:=F_i^n/[n]!_{d_i},\ K_i^{\pm 1},\
i=1,\ldots ,r,\ n\ge 0.
$$
It is analogous to the classical
integral form of the universal enveloping algebras
due to Kostant.

By definition set $U_3=k\otimes_{\Q [v,v^{-1}]}U_{\Z}$. Then $U_3$
contains elements $E_i ,F_i,K_i$
and $K_i ^{-1},\ i=1,\ldots,r,$ that
satisfy the basic relations for the generators
of $U_2$. Thus
there exists a homomorphism $f: U _2\map  U_3$
mapping the generators of
$U_2$ to the corresponding elements of
$U_3$.

Its image is a finite dimensional subalgebra
$\u \subset U_3$ (see   [L1],
8.2). The kernel of $f$ contains $I$, and we obtain the surjective map:
$U\map \u$.  Both algebras are
$X$-graded, and the map preserves the grading.

\subsubsection{}
\Lemma (see   [AJS],
1.3)
 \label{decomp}
 The following statements hold for $U$ and $\u$:\\
\qquad(i) $U=U^{+}\otimes U^0\otimes U^{-}$ as a vector space; \\
\qquad(ii) $\u=\u^{+}\otimes \u^0\otimes \u^-$ as a vector space;\\
\qquad(iii) both maps $U^{\pm }\map \u^{\pm }$ are
isomorphisms of algebras;\\
\qquad(iv) $\u^0$ is equal to
the quotient algebra  of $U^0$ by the ideal generated by all
$K_i^{2p}-1$.$\qed $
\ms
Thus we obtain PBW-type bases in
$\u$. We call both
$U$ and $\u$ the finite quantum groups. The algebra $\u^-\otimes\u^0$ (resp.
$\u^0\otimes\u^+$) is called the negative
(resp. the positive) Borel subalgebra in $\u$
and is denoted by $\b^-$ (resp. by $\b^+$).

\subsection{}
Using PBW-type bases, Kac, De Concini and
Procesi defined some remarkable filtrations on $U_2$,
$U$ and $\u$ that generalize the usual
PBW-filtrations on universal enveloping algebras (see  [DCKP]).

Consider the lexicographically ordered set
$S=\Z_{+}^{2N+1}$, where $N$ is the number of the positive
roots.  A filtration of a vector space by subspaces numbered by the
set $S$ is called $S$-filtration.

\subsubsection{}
We fix a convex order on the set of the positive roots
$R^{+}$. Roughly speaking the convex
property means that the $q$-commutator of two root vectors
$E_{\alpha}$ and $E_{\beta}$ in $U_1$ consists
of monomials formed only from root vectors between $\alpha$ and $\beta$
in the order (see e.g. [DCKP] for the exact definition).  We denote the
monomial
$$
E_{\beta_{1}}^{s_{1}}\ldots E_{\beta_{N}}^{s_{N}}K_1^{b_1}\ldots K_r^{b_r}
F_{\beta_{1}}^{t_{1}}\ldots F_{\beta_{N}}^{t_{N}}
$$
by
$M_{(s,t,\alpha)}$, where $\alpha =\sum_{j=1}^r b_j\alpha_j$.
We define the total height
$$
d_{0}(M_{(s,t,\alpha)}) =
\sum_{i=1}^N (s_i+t_i)\hgt\beta_i,
$$ where
$\hgt\beta$ is the height of the root $\beta$,
and the total degree
$$
d(M_{(s,t,\alpha)})=(s_N,\ldots,s_1,t_1,\ldots,t_N,d_0(M_{(s,t,\alpha)}))
\in\Z_{+}^{2N+1}.
$$
We introduce the $S$-filtration $F$ on $U_2$ by the total degree.

\subsubsection{}
\Lemma (see e.g. [DCKP], 4.2)\label{rel}
The associated polygraded algebra $\gr U_2$ for the $S$-filtration on
$U_2$ is generated by the elements $E_\alpha ,F_{-\alpha}$
($ \alpha \in R^+$) and
$K_i^{\pm 1}$ ($i=1,\ldots ,r$) satisfying the following relations:
$$
K_iE_\beta =\zeta ^{(\alpha_i|\beta )}E_\beta K_i,\
K_iF_{-\beta}=\zeta ^{-(\alpha_i|\beta )}F_{-\beta}K_i,\ \beta\in R^+;
$$
$$
K_iK_i^{-1} =K_i^{-1}K_i=1,\ K_iK_j=K_jK_i;
$$
$$
E_\alpha F_{-\beta}=F_{-\beta }E_\alpha ,\ \alpha ,\beta \in R^{+};
$$
$$
E_\alpha
E_\beta =\zeta ^{-( \alpha | \beta )}E_\beta E_\alpha,\
F_{-\alpha }F_{-\beta }=\zeta ^{( \alpha | \beta )}F_{-\beta}F_{-\alpha }
$$
if $\alpha ,\beta \in R^{+}$ and $\alpha >\beta $ in the
convex order on $R^{+}.\qed$
\bs
The filtration $F$ defines
$\Z_{+}^{2N+1}$-filtrations on $U$ and $\u$ and
$\Z_{+}^{N+1}$-filtrations on $\u^{\pm}$, $\u^{+}\otimes \u^{0},\
\u^{0}\otimes \u^{-}$ etc.  We denote them by $F$ as well.

\subsubsection{}
\Cor \label{fullrel}
The graded algebra $\gr^{F}U$ is
generated by the elements $E_\alpha$, $F_{-\alpha}$, $\alpha\in R^+$, and
$K_i^{\pm 1}$, $i=1,\ldots,r$,
subject to the relations from  Lemma~\ref{rel}  and the following
relation:  $E_\alpha ^p=F_{-\alpha}^p=0\text{ for every }
\alpha \in R^+.\qed $
\bs
Recall that a finite dimensional associative algebra $A$ is called
Frobenius if the dual module to the right regular module over $A$ is
isomorphic to the left regular $A$-module:  $A_L\cong \Hom_k(A_R,k)$.
The nondegenerate     bilinear pairing $A\times A\map k$ induced by
this identification is called the trace on $A$.

\subsubsection{}
\Lemma \label{frob}
Let $A$ be  a finite dimensional filtered algebra such that the top
component of the  corresponding graded algebra is one dimensional, and
$\gr A$ is Frobenius with the trace defined as follows:
$$
(x,y)\mapsto (x\cdot y)_{top}\in (\gr A)_{top}\cong k.
$$
Here
$\cdot$ denotes  the multiplication in $\gr A$. Then $A$ is also
Frobenius.  The trace on $A$ is defined as follows:  $$(x,y)\mapsto
xy\mapsto (xy)_{top}\in (grA)_{top}\cong k.$$ The product here is the
product in $A$.  $\qed$

\subsubsection{}
\Lemma
The algebras $\u,\ \u^+$ and $\u^-$
are Frobenius.
\ms
\dok
By \ref{fullrel} the algebras $\gr \u^{\pm}$
satisfy the conditions of the Lemma~\ref{frob}. Thus $\u^{\pm}$ are Frobenius.
The trace on $\u$ is defined in
[Xi], 2.9.$\qed $

\subsubsection{}
We will need a filtration on $\u$ that differs a little from $F$.
We set
$$
d_1(M_{(s,t,\alpha)})=\sum_{j=1}^N s_j\hgt\beta_j
$$
and  introduce the partial degree
$$
\widetilde{d}(M_{(s,t,\alpha)})=(s_N,\ldots,s_1,
d_1(M_{(s,t,\alpha)}))
$$.
Let $F'$ be the $\Z^{N+1}$-filtration of $\u$ by the partial
degree. Note that $F'$ coincides with $F$ on $\u^+$ and it coincides on $\u^-$
with the natural $\Z$-grading obtained from the $Y^-$-grading. Thus
$\gr^{F'}(\u^-)=\u^-,\ \gr^{F'}(\u^+)=\gr^F(\u^+).$

Recall that an augmented subalgebra $B\subset A$ is called normal if
the left and the right ideals in $A$ generated by the augmentation
ideal $\overline{B}$ in $B$ coincide. Then $A//B$ denotes the quotient
algebra of $A$ by the two sided ideal. The algebras $\u,\ \u^{\pm}, \
\gr^F(\u )$ and $\gr^F(\u^{\pm})$ are naturally augmented: the
augmentation is provided by the map
$$
E_\alpha\mapsto 0,\
F_{-\alpha}\mapsto 0,K_i\mapsto 1,\ \alpha\in R^+,\ i=1,\ldots ,r.
$$

\subsubsection{}
\Lemma \label{filtration}
$F'$ defines a filtration
on $\u$ compatible with the multiplication in $\u$.
$$
\gr^{F'}(\u )= \u^{-}\otimes \u^{0}\otimes \gr^{F}(\u^{+})
$$
as a vector space, $\u^{-},\ \u^0$ and $\gr^{F}(\u^+)$ are
subalgebras in $\gr^{F'}(\u )$. Finally, $\gr^F(\u^+)$
is normal in $\gr^{F'}(\u )$ and elements of $\u^-$ commute with elements of
$\gr^{F'}(\u^+)$.
\ms
\dok
The first statement is an easy consequence of  Lemma~\ref{rel}  since
if  two monomials $M_{(s,t,\alpha)}$ and $M_{(u,v,\beta)}$ are of equal total
degrees then their partial degrees also coincide:
$\widetilde{d}(M_{(s,t,\alpha)})= \widetilde{d}(M_{(u,v,\beta)})$.

As $F$ and $F'$ coincide on
$\u^+$, $\gr^F(\u^{+})$ is a subalgebra in
$\gr^{F'}(\u )$.

The decomposition of $\gr^{F'}(\u )$
into tensor product follows from the similar statement for $\u$
(see   Lemma~\ref{decomp}).

To prove that elements of $\u^-$
commute with elements of $\gr^F(\u^+)$ in $\gr^{F'}(\u)$ note that for each
$i=1,\ldots ,r$
$$
[ F_i,M_{(s,0,0)} ]\in \u^0\otimes \u^+ \text{ and }
d_1([F_i,M_{(s,0,0)}])<d_1(M_{(s,0,0)}).\qed
$$

\subsubsection{}
\Lemma
The algebra
$\gr^{F'}(\u )$ is Frobenius.
\ms
\dok
As a vector space $\gr^{F'}({\u})=\u^-\otimes \u^0\otimes \gr^F(\u^+).$
Set the linear form, inducing the
trace on $\gr^{F'}({\u})$, equal to the tensor product of the
linear forms inducing the traces on the componenents.  $\qed$

\subsection{}
Now we are going to describe the categories of $\u$-modules we will work with.
$\C$ is the category of
$X$-graded left $\u$-modules $M=\underset {\lambda \in
X}{\bigoplus}\text{ }_{\lambda}M$,\text{ }$\dim \text{
}_{\lambda}M<\infty$, such that $K_i$ acts on $ _{\lambda}M$ by
multiplication by $\zeta^{(\alpha_i^{\vee}|\lambda)}$,
where $\alpha_i^{\vee}$ are the simple coroots,
$$
E_i:\ _{\lambda}M\map \ _{\lambda +\alpha_i}M,\
F_i:\ _{\lambda}M\map \ _{\lambda -\alpha_i}M.
$$
Morphisms in
$\C$ are the morphisms of $X$-graded $\u$-modules that
preserve gradings.  $\C^{(r)}$ is the category of  right
$X$-graded    $\u$-modules $N=\underset {\lambda \in
X}{\bigoplus}\text{ }_{\lambda}N$, $\dim \text{ }_{\lambda}N<\infty$, such
that $K_i$ acts on $_{\lambda}N$ by multiplication by $\zeta^{-(
\alpha_i^{\vee}|\lambda )}$,
$$
E_i:\ _{\lambda}N \map \ _{\lambda+\alpha_i}N,\
F_i:\ _{\lambda}N\map  \ _{\lambda-\alpha_i}N,
$$
with morphisms that preserve $X$-gradings. One can define the categories
$\C (\b )$, $\C (\gr^{F'}(\u ))$ and $\C^{(r)}(gr^{F'}(\u))$ in a similar way.

We define the twisting functors  by
elements of the weight lattice  $\lambda \in pX$ on the category $\C$:
$$
M\mapsto M\langle\lambda\rangle\text{, where } _{\mu}M\langle\lambda\rangle:=
\ _{\lambda+\mu}M,
$$ with the same action of $\u$.
We denote the space $\underset{\beta \in X}{\bigoplus }
\Hom_{\cal C}(L,M(p\beta ))$  by $\Hom_{\u}(L,M)$.

We will also need the categories of finite dimensional left
$\u^{\pm}$-modules denoted by ${\u}^{\pm}-\Mod$.

We call the
${\u}$-module $M_{\lambda}^+:=\u\otimes_{\b^-}k_{\lambda}$ (resp. the
$\u$-module $M_{\lambda}^-:=\u\otimes_{\b^+}k_{\lambda}$)
the positive (resp. negative) left Verma module with the highest (resp.
lowest) weight $\lambda$. Here $k_{\lambda}$ is one-dimensional
$\u^0$-module placed in the
$X$-grading $\lambda$, the trivial action of ${\u}^{\pm}$ equips it with the
structure of a $\b^{\pm}$-module.

Contragradient Verma modules $M^{- *}_{\lambda}$ are defined in the
following way:  $M^{- *}_{\lambda}=\Hom_{\b^-}({\u},k_{\lambda})$
with the natural left action of ${\u}$.
We will need the following statement.

\subsubsection{}
\label{verma}
\Lemma (see  [AJS], 4.10, 4.12)

\qquad(i) $M_{\lambda}^{+}=M_{\lambda+(p-1)2\rho}^{- *}$;

\qquad(ii) $\Hom_{\C}(M_{\lambda}^-,M_{\mu}^+)=0$ if $\lambda \ne
\mu+(p-1)2\rho$;

\qquad(iii) $\Hom_{\C}(M_{\lambda}^-,M_{\lambda+(p-1)2\rho}^+)=k$.
$\qed$

\subsubsection{}
Every Verma module
$M_{\lambda}^-$ has a unique simple  quotent module $L_{\lambda}$ with
the highest weight $\lambda \in X$, and this way one obtains the full
list of simple modules in the category $\C$ (see  [AJS], 4.1).

\section{Semiinfinite cohomology of modules over associative algebras}

Consider a free abelian group $X$ of the finite rank $r$ and  its
subgroup $Y$ generated by a set $\Sigma \subset X$ consisting of $r$
elements, such that the elements of $\Sigma $ form a base of the
vector space $X\otimes \Q$.  Denote the subsemigroup in $Y$ generated
by the set $\Sigma$ (resp. by the set $-\Sigma$)  by $Y^+$ (resp. by
$Y^-$). In this section we do not suppose that $X$ and $Y$ are the
weight and the root lattices corresponding to some root system, but
our notation is adopted to that case.

Let $v:\ X\map\Q$ be a linear function defined as follows: for every
$\alpha\in\Sigma$ $v(\alpha)=1$, $v$ is extended on $X$ by linearity.

\subsection{}
\label{setup}
Suppose we have an $Y$-graded associative
algebra $A$ with the $Y$-graded subalgebras $B$ and $N$ satisfying
 the following conditions :\\
\qquad(i) $N$ is graded by $Y^{+}$;\\
\qquad(ii) $N_0=k;$\\
\qquad(iii) $\dim N_\beta<\infty$ for any $\beta \in
Y^{+};$\\
\qquad(iv) $B$ is graded by $Y^{-}$;\\
\qquad(v) the multiplication in $A$
defines the isomorphisms of $Y$-graded vector spaces
$$
B\otimes N\map A\text{ and }N\otimes B\map A.
$$
In particular $N$ is naturally augmented.

Note that our conditions are satisfied both for finite quantum groups and for
universal enveloping aslgebras of graded Lie algebras. In the latter case
the decomposition into tensor product of negative and positive parts is
given by the natural decomposition of the Lie algebra into the direct sum
of its negatively and positively graded parts.

\subsection{}
Consider the category $A\mod$ of $X$-graded left $A$-modules with
morphisms that preserve $X$-gradings. The corresponding category of right
$A$-modules is denoted by $\mod -A$. We will need the following subcategories
in the category of complexes $ {\cal K} om(A\mod )$.

For an $X$-graded module $M$ denote by $\supp M$ the  set
$\{\alpha \in X\mid \text{ }_{\alpha}M\neq 0\}$. Denote by $X_\Q^{+}(\beta)$
(resp. $X_\Q^{-}(\beta)$) the convex cone in $X\otimes \Q$,
generated by $\Sigma $ (resp.  $-\Sigma$) with the vertex in $\beta \in
X\otimes k$.

For $s_1,s_2,t_1,t_2\in\Z$, $s_1,s_2>0$, the set
$$
 \{(p,q)\in\Q^{\oplus 2}|s_1p+q\ge t_1,\ s_2p-q\ge t_2\}
$$
(resp. the set
$$
 \{(p,q)\in\Q^{\oplus 2}|s_1p+q\le t_1,\ s_2p-q\le t_2\})
$$
is denoted by $\Q^{\uparrow}(s_1,s_2,t_1,t_2)$ (resp. by
$\Q^{\downarrow}(s_1,s_2,t_1,t_2)$).

\label{wedge}
Consider the category $\upC (A)$
(resp.  $\dC (A)$) in ${\cal K} om(A\mod)$,
consisting of complexes of $X$-graded
$A$-modules $M=\underset{q\in \Z,\lambda \in X}{\bigoplus }
\text{ }_{\lambda}M^q$ satisfying the following conditions:

(U) there exist  $s_1,s_2,t_1,t_2\in\Z,\ s_1,s_2>0,$ such that
$v(\supp M\bul)\subset\Q^{\uparrow}(s_1,s_2,t_1,t_2)$ and for any
$(p,q)\in v(\supp M\bul)$ the set $v^{-1}(p,q)$ is finite;
resp.

(D) there exist  $s_1,s_2,t_1,t_2\in\Z,\ s_1,s_2>0,$ such that
$v(\supp M\bul)\subset\Q^{\downarrow}(s_1,s_2,t_1,t_2)$ and for any
$(p,q)\in v(\supp M\bul)$ the set $v^{-1}(p,q)$ is finite.

Here the set $v(\supp M\bul)$ is considered as a subset of a $(p,q)$-plane:
$v(\lambda,q):=(v(\lambda),q)$.

\subsection{}
Now we are going to define the standard complex for the computation of
semiinfinite cohomology of $A$-modules. Consider the right $N$ module
$$
N^*=\Hom_k(N,k)=\underset {\beta \in X^+}{\bigoplus}\Hom_k(\ _{\beta}N,k).
$$
The right action of $N$ is defined as follows:
$$
n\cdot f(m)=f(nm),\ n,m \in N,\ f\in N^*.
$$
Voronov calls the right $A$-module $S_A=N^*\otimes _NA$ the right semiregular
representation (see [V]). Obviously $S_A$ is isomorphic  to $N^*\otimes B$ as a
right $B$-module.

For two $X$-graded $B$-modules $L$ and $M$ denote the space
$\underset{\lambda \in
X}{\bigoplus}\Hom_{B\mod}(L,M\langle\lambda\rangle)$, where
$\langle\lambda\rangle$ is the grading shift, by $\Hom_{B}(L,M)$.

\subsubsection{}
\Lemma $S_A=\Hom_{B}(A,B)$ as a right $A$-module.
\ms
\dok Define the pairing $\phi :\ S_A\times A\map B$ as follows:
$$
\phi(f\otimes a_1,a_2)=f_1(a_1a_2),
$$
where $f_1$ denotes $f$ used by the first argument in $A\otimes B$.
The required isomorphism is provided by $\phi$.$\qed$

Set $A^?=\operatorname{End}_A(S_A)$. The functors of induction and coinduction
provide the natural inclusions of algebras $N\hookrightarrow A^?$ and
$B\hookrightarrow A^?$. Clearly up to a certain completion $A^?=B\otimes N$
as a vector space.

\subsubsection{}
\Lemma The subspace $\oppA=B\otimes N\subset A^?$ is a $Y$-graded subalgebra.
$\qed$
\bs
Thus $S_A$ becomes an $\oppA -A$ bimodule.

\subsubsection{}
Recall the bar construction for the algebra $A$ with respect to
the subalgebra $B\subset A$.
The standard bar resolution $\tilBarb (A,B,M)\in {\cal K}om(A\mod)$
of an $A$-module $M$ is defined as follows:
$$
\tilBar^{-n}(A,B,M)=A\otimes_B\ldots \otimes_BA\otimes_BM \
(n+1 \text{ times}),
$$
$$
\operatorname{d}(a_0\otimes \ldots \otimes a_n\otimes v)=
\sum_{s=0}^{n-1}(-1)^s a_0
\otimes\ldots\otimes a_sa_{s+1}\otimes\ldots\otimes v
+(-1)^na_0\otimes\ldots\otimes a_{n-1}\otimes a_nv.
$$
Here $a_0,\ldots ,a_n\in A,\ v\in M$.

\subsubsection{}
\Lemma The subspace $\linBarb (A,B,M):$
$$
(\linBar)^{-n}(A,B,M)=\{ a_0\otimes\ldots\otimes a_n\otimes v\in
\tilBar^{-n}(A,B,M)|\ \exists \ s\in\{1,\ldots,n\}:\ a_s\in B\}
$$
is a subcomplex in $\tilBarb (A,B,M)$.$\qed$
\ms
The quotient $\Barb (A,B,M)=\tilBarb (A,B,M)/\linBarb (A,B,M)$ is called the
restricted bar resolution of the $A$-module $M$ with respect to the subalgebra
$B$.
The ideal of augmentation in $N$ is denoted by $\overline{N}$.

\subsubsection{}
\Lemma

\qquad(i) $\Bar^{-n}(A,B,M)=N\otimes \overline{N}^{\otimes n}\otimes M$
as a left $N$-module;

\qquad(ii) $\Barb(A,B,M)=\Barb(N,k,M)$ as a complex of $N$-modules. In
partiqular $\Barb (A,B,M)$ is a $N$-free resolution of the $A$-module $M$.
$\qed$
\bs
For any $M\bul\in \upC (A)$ consider the total complex of its bar resolution
$\Barb (A,B,M\bul)$ and the complex of $\oppA$-modules
$S_A\otimes_A\Barb (A,B,M\bul )$
\bs
{\bf Remark:\ } Since $S_A$ is $B$-free and $\Barb (A,B,M\bul )$ is $N$-free
$$
H\bul (S_A\otimes_A\Barb (A,B,M\bul ))=\Tor^A(S_A,M\bul ).
$$
For a complex of left $\oppA$-modules $L\bul\in\dC (\oppA )$ we denote
the total complex of its restricted bar resolution with respect to the
subalgebra $N\in\oppA$ by $\Barb (\oppA ,N,L\bul )$.

\subsubsection{}
\Def Let $L\bul \in\dC (\oppA ),\ M\bul \in \upC(A)$. The standard
complex for the computation of semiinfinite $\ext$ functor $\stand
(L\bul ,M\bul )$ is defined as follows:
$$
\stand (L\bul ,M\bul )=
\Hom\bul_{\oppA}(\Barb (\oppA ,N,L\bul ),S_A\otimes_A\Barb (A,B,M\bul )).
$$
By definition set $\Exts_A(L\bul ,M\bul )=H\bul (\stand (L\bul
,M\bul ))$.

Note that unlike usual $\ext$ and $\tor$ functors semiinfinite cohomology
exists both in negative and positive homological degrees even for $L$ and $M$
being complexes-objects (see [GeM]).

As a vector space
$\stand (L\bul ,M\bul )=\underset{n,m}{\bigoplus}
\Hom_k(\overline{B}^{\otimes n}\otimes L\bul ,
\overline{N}^{\otimes m}\otimes M\bul )$. Here $\overline{B}$ denotes the
space $B/k$.
Since both arguments of the semiinfinite $\ext$ functor are $X$-graded,
both the standard complex $\stand (L\bul ,M\bul )$ and its cohomology are
also $X$-graded:
$$
\ _{\gamma}\stand (L\bul ,M\bul )=
\underset{\alpha-\beta=\gamma;n,m}{\bigoplus}\Hom_k(
\ _{\beta}(\overline{B}^{\otimes n}\otimes L\bul ),
\ _{\alpha}(\overline{N}^{\otimes m}\otimes M\bul )),
$$
$$
\ _{\gamma}\Exts_A(L\bul ,M\bul )=H\bul (\ {\gamma}\stand (L\bul ,M\bul )).
$$

From now on the zeroth $X$-grading component of $\Exts_A(L\bul ,M\bul $) is
denoted by $\Exts_{A\mod }(L\bul ,M\bul )$.  \label{0}

\subsection{}
Consider the filtrations $^{(I)}F$ (resp.
$^{(II)}F$) on $\stand (L,M),\ L\in \oppA \mod ,M\in A\mod ,$
by the number $n$ (resp. $m$).
The $E_0$-terms  of the corresponding spectral sequences are as follows:
\begin{multline*}
^{(I)}E_0^{p,q}=\Hom_{\oppA}(
\oppA\otimes_N\underbrace{(\oppA /N)\otimes_N\ldots\otimes_N(\oppA /N)}_{p}
\otimes_NL,S_A\otimes_A\Bar^{q}(A,B,M)) \\=
\Hom_N(\underbrace{(\oppA /N)
\otimes_N\ldots\otimes_N(\oppA /N)}_{p}\otimes_NL,
S_A\otimes_B
\underbrace{(A/B)\otimes_B\ldots\otimes_B(A/B)}_{-q}\otimes_BM) \\=
\Hom_k(\overline{B}^{\otimes p}
\otimes L,\underbrace{(A/B)\otimes_B\ldots\otimes_B
(A/B)}_{-q}\otimes_BM) \\=
\Hom_B(\Bar^{-p}(B,k,L),\underbrace{(A/B)\otimes_B\ldots\otimes_B(A/B)}_{-q}
\otimes_BM)
\end{multline*}
with the differential being that in $\Barb (B,k,L)$, and
$$
^{(II)}E_0^{p,q}=\Hom_k(\underbrace{(\oppA /N)\otimes_N\ldots
\otimes_N(\oppA /N)}_{p}\otimes_NL,k)\otimes_N\Bar^{-q}(N,k,M)
$$
with the differential being that in $\Barb (N,k,M)$.
Thus the $E_1$ terms are as follows:
\begin{gather*}
^{(I)}E_1^{p,q}=\ext^p_B(L,
\underbrace{(A/B)\otimes_B\ldots\otimes_B(A/B)}_{-q}\otimes_BM),\\
^{(II)}E_1^{p,q}=\tor^N_{-q}(\Hom_k(
\underbrace{(\oppA /N)\otimes_N\ldots\otimes_N(\oppA /N)}_{p}\otimes_NL,k)
,M).
\end{gather*}
Both spectral sequences are $X$-graded.

\subsubsection{}
\Lemma
Let $L\in \oppA- \mod,\ \supp L\in X_\Q^-(\lambda),\
M\in A\mod,\ \supp M\in X^+_\Q(\mu )$. Then for a  fixed $\beta \in
X$ both $\ _{\beta}(^{(I)}E^{p,q})$ and $\ _{\beta}(^{(II)}E^{p,q})$
converge.
\ms
\dok
Since $B$ is $Y^-$-graded, there exists $\beta_0\in X$ such that for every
$p\ge 0\ \supp \overline{B}^{\otimes p}\otimes L$ belongs to $X^-_\Q(\beta_0)$.
The $X$-graded space $\underset{q\ge 0}{\bigoplus}
\overline{N}^{\otimes q}\otimes M$
satisfies the condition (U). Thus in a fixed $X$-grading component $\beta$
both spectral sequences of the complex $\ _{\beta}\stand (L,M)$ are situated in the part of the
$(p,q)$-plane wich is bounded in $p$ from the left and in $q$
both from the left and from the right.
$\qed$

Recall that an object $M\in A\mod$ is injective (resp. projective)
relative to the subalgebra
$N$ if for every complex of $A$-modules $C\bul$ such that $C\bul$ is homotopic
to zero as a complex of $N$-modules $H\bul (\Hom\bul_A(C\bul ,M))=0$.
(resp. $H\bul (\Hom\bul_A(M,C\bul )=0$).

\subsubsection{}
\Lemma \label{main} The following facts hold
for $L\bul\in \dC (\oppA ),$ $M\bul\in
\upC (A)$:

\qquad(i) if $M\bul$ is $N$-projective, then
$$
\Exts_A(L\bul ,M\bul )=H\bul (\Hom_{\oppA }(\Barb(\oppA ,N,L\bul ),
S_A\otimes_AM\bul));
$$
\qquad(ii) if $L\bul$ is $B$-projective, then
$$
\Exts_A(L\bul ,M\bul )=H\bul (\Hom\bul_{\oppA }(L\bul ,S_A\otimes_A
\Barb (A,B,M\bul )));
$$
\qquad(iii) if $M\bul$ is both $N$-projective and
$A$-injective relative to $N$,
then
$$
\Exts_A(L\bul ,M\bul )=
H\bul (\Hom_{\oppA }\bul (L\bul ,S_A\otimes_AM\bul )).
$$
\ms
\dok
(i) Consider the canonical mapping $\varphi:\ \Barb (A,B,M\bul )\map M\bul$.
Then $\Cone\bul \varphi$ is an exact complex of $N$-projective $A$-modules
satisfying (U). In partiqular it is homotopic
to zero as a complex of $N$-modules.
As a complex of $N$-modules
$$
S_A\otimes_A\Cone\bul\varphi=N^*\otimes_N\Cone\bul\varphi,
$$
where $N^*$ is considered as a $N$-bimodule. Thus
$S_A\otimes_A\Cone\bul\varphi$ is also homotopic to zero over $N$. Now
by Shapiro lemma $\oppA$-modules induced from $N$-modules are
relatively projective, so $\Barb (\oppA ,N,L\bul)$ consists of
relatively projective modules. Let
$$
 \widetilde{\varphi}:\ \stand
(L\bul ,M\bul )\map \Hom\bul_{\oppA }( \Barb (\oppA ,N,L\bul
),S_A\otimes_AM\bul )
$$
be the  morphism complexes corresponding to
$\varphi$. Then
$$
 \Cone\bul\widetilde{\varphi}=\Hom\bul_{\oppA}(\Barb
(\oppA ,N,L\bul), S_A\otimes_A\Cone\bul\varphi ).
$$
We prove that
the latter complex is exact. Consider the bigrading on
$\Cone\bul\widetilde{\varphi}$:
$$
\Cone^{m,n}\widetilde{\varphi}=\underset{p+s=m}{\bigoplus} \Hom_{\oppA }
(\Bar^p(\oppA ,N,L^s),\Cone^n\varphi).
$$
Our grading conditions
provide that the spectral sequence of the bigraded complex converges.
On the other hand
$$
E^{m,\bullet}_1=\underset{p+s=m}{\bigoplus}
H\bul(\Hom\bul_{\oppA }(\Bar^p(\oppA ,N,L^s),\Cone\bul\varphi ))=0.
$$
(ii) The proof is similar to the previous one.

(iii) It is sufficient to prove that the mapping
$$
\Hom\bul_{\oppA}(L\bul,S_A\otimes_AM\bul )\map
\Hom\bul_{\oppA}(\Barb (\oppA ,N,L\bul ),S_A\otimes_AM\bul )
$$
is a quasiisomorphism.
We are going to show that if $M$ is both $N$-projectine and injective
relatively to $N$ then $S_A\otimes_AN$ is injective over $\oppA$.

First note that the functor $S_A\otimes_A*$ takes $N$-free
modules to $N$-cofree ones,
thus it takes $N$-projectives to $N$-injectives.

The functor $\Hom_{\oppA}(S_A,*)$ is the right conjugate functor
for $S_A\otimes_A*$. It is left exact and is
well defined on $N$-modules since it can be written as follows:
$M\mapsto \Hom_N(N^*,M).$ Thus $S_A\otimes_A*$ preserves
relative injectiveness.

Finally note that a $\oppA$ module that is both
$N$-injective  and relatively injective is also $\oppA$-injective.
Thus for a finite exact complex $P\bul$  the complex
 $\Hom\bul_{\oppA}(P\bul ,S_A\otimes_AM\bul )$ is exact. It remains to
check the convergence of a spectral sequence similar to the one
from the first statement.
$\qed$

\subsection{} \label{derived}
\Theorem Semiinfinite $\ext$ functor is well defined on the
corresponding derived categories:
$$
\Exts_A:\ {\cal
D}^{\downarrow}(\oppA )\times {\cal D}^{\uparrow}(A)\map {\cal D}
({\cal V}ect).
$$
Here ${\cal D}^{\downarrow}(\oppA )$ (resp. ${\cal
D}^{\uparrow}(A)$) denotes the localization of the category $\dC
(\oppA )$ (resp. $\upC (A)$) by the class of quasiisomorphisms.
\ms
\dok
We are to prove that $\Exts_A(L\bul ,M\bul )=0$ for
$L\bul\in\dC(\oppA ),\ M\bul \in \upC (A)$ if either of the arguments
is exact. Suppose $M\bul $ is exact, the proof in the other case is
quite similar.

We fix $\lambda\in X$.
Consider the following bigrading on $\stand (L\bul ,M\bul )$:
$$
C^{\frac \infty 2 \ p,q}(L\bul ,M\bul ) = \underset{m-n+s=p}{\bigoplus}
\Hom_{k}((\overline{B})^{\otimes m}\otimes L^{-s}),
(\overline{N})^{\otimes n}\otimes M^q).
$$
The spectral sequence of the bicomplex
$\ _{\lambda}C^{\frac\infty2\bullet\bullet}(L\bul,M\bul)$ converges since
it is situated in the  part of the $(p,q)$-plane bounded in $p$ from
the left and in $q$ both from the left and from the right. Thus the
total complex $\stand (L\bul, M\bul )$ is exact.  $\qed$

\subsubsection{}
{\bf Remark:} The last two statements
show that one can use an arbitary resolution
of $M\bul$ that is both $N$-projective and
injective relatively to $N$ for the computation
of $\Ext_A(L\bul ,M\bul )$.

\section{Calculation of the character of
$\protect\Exts_{\C}(k,k)$.}

In this section we are going to define a $\g$-module structure  on
the semiinfinite cohomology of the trivial $\u$-module. We will
also calculate the character of this $\g$-module.

From now on $X$ and $Y$ denote the weight and the root lattice respectively,
the linear function $v$ coincides with $\hgt$.

We will need the
subcategories  in the category of complexes ${\cal K}{om}(\C )$,
satisfying the condition (U) (resp. (D)) from the previous section
(the category ${\C}$ is defined in 2.4).  These categories are denoted
by $\upC$ (resp. by $\dC$).

\subsection{}
Fix the triangular decomposition of the finite quantum
group:  $\u=\b^-\otimes\u^+$, i.~e. $A=u$, $B=\b^-$, $N=\u^+$ in the
notations of~\ref{setup}.

\subsubsection{}
Lemma: \label{dual}
$\u^?=\u^{\sharp }=\u.$
\ms
\dok
The algebra
$\u^+$ is Frobenius (see   \ref{frob}), thus the right
$\u$-module $S_{\u}$ is isomorphic to the
right regular $\u$-module.  But for any finite dimensional
algebra the algebra of endomorphisms of the right regular module is
isomorphic to the algebra itself.$\qed$

Note that the category $\C$ differs from the category $\u\mod$ thus to
define semiinfinite cohomology of $\u$-modules one has either to introduce
a $X$-graded algebra $A$ such that $\C\cong A\mod$ (the  algebra  $A$ is
constructed in particular in [AJS],Remark 1.4) and to consider
$\Exts_{A\mod}(*,*)$ or to define semiinfinite cohomology in the category
$\C$ explicitly.

Note that for $M\in\C$ the complexes $\Barb(\u,\b^+,M)$ and $\Barb(\u,\b^-,M)$
belong to $\cal{K}om(\C)$.

\subsubsection{}
\Def
For $L\bul\in\dC,\ M\bul\in\upC$ the semiinfinite $\ext$ functor is
defined as follows:
$$
\Exts(L\bul, M\bul):=
H\bul\left(\Hom\bul_{\C}(\Barb(\u,\b^+,L\bul),\Barb(\u,\b^-,M\bul))\right).
$$
Note that by Lemma~\ref{dual} the module $S_{\u}$ is isomorphic to
the right regular $\u$-module, so the definition is a direct analogue
of~\ref{0}. In partiqular the statements of Lemma~\ref{main} remain
true.

\subsubsection{}
The subalgebra $\u^+\subset\b^+$ (resp. $\u^-\subset\b^-$) is normal.
Clearly $\b^+//\u^+=\u^0$ is semisimple being the group algebra
of the group $\Z /2p\Z$. In particular a $\u$-module $L$ is $\b^+$-projective
(resp. $\b^-$-projective) if and only if it is $\u^+$-projective
(resp. $\u^-$-projective).

\subsection{}
 Consider the following  $\u^+$-free $\u^-$-cofree resolution of a $U_3$-module
 $M\bul\in\upC (U_3)$:
$$
R\bul (M\bul ):=\Barb (U_3,B_3^+,k)^*\otimes \Barb (U_3,B_3^-,M\bul ).
$$
Here $B_3^{\pm}\subset U_3$ denotes the positive and the negative
Borel subalgebras in the "big" quantum group $U_3$. The definition of
the tensor product over the base field uses the standard Hopf algebra
structure on $U_3$. The left $U_3$-module structure on
$$
\Barb
(U_3,B^+_3,k)^*:=\underset{\lambda \in X} {\bigoplus}\Hom_k(\
_{\lambda}\Barb (U_3,B_3^+,k), k)
$$
is defined using the antipode in
$U_3$.

Evidently $R\bul (M\bul )$ satisfies the condition of Lemma~\ref{main}(iii),
in particular it satisfies the condition (U).
Thus it can be used for the computation of semiinfinite cohomology of
$\u$-modules.

\subsubsection{}
\Lemma Let $M\bul\in \upC (U_3)$. Then
$\underset {\lambda\in X}{\bigoplus}\Exts_{\C}(k,M\bul \langle
p\lambda\rangle)$ admits a structure of  $\frak g$-module.
\ms
\dok
In [L1], Theorem 8.10, it is proved that the algebra ${\u}$ is normal
in $U_3$, and the quotient algebra $U_3//{\u}$ is isomorphic to the
 universal enveloping algebra of the semisimple Lie algebra $\g$.
Thus by Shapiro lemma
$$
\underset {\lambda \in X}{\bigoplus}
\Exts_{\C}(k,M\bul \langle p\lambda\rangle)=H\bul (\Hom_{\u}\bul
(k,R\bul (M\bul)))= H\bul(\Hom_{U_3}\bul (U_3//\u ,R\bul (M\bul ))).
$$
The left $U_3$-module $U_3//{\u}$  is naturally equipped with a
right action of the quotient algebra $U_3//\u=U(\g )$ commuting with
the left action of $U_3$. Thus the semiinfinite $\ext$ spaces carry
the natural structure of $\g$-modules.  ${\qed }$

\subsection{}
The rest of this section is devoted to the computation of the
character  of $\Exts_{\C}(k,k\langle p\lambda \rangle)$.
The  problem  is that we do not know ``the minimal" $\u^-$-free
resolution of the trivial ${\u}$-module $k$.

\subsubsection{}
{\bf Conjecture:} There exists  a resolution
$R_{\min }^{ -\bullet}(k)$ of the trivial ${\u}$-module $k$,
filtered by Verma modules $M_\lambda ^{-}$, such that the character of the
space spanned by the  highest weight vectors $v_{\lambda}\in
M^-_{\lambda}$ is given by the formal series
$$
ch(t):=\frac{\underset{w\in
W}{\sum }e^{w(\rho )-\rho }t^{l(w)}}{\underset{\alpha \in R^{+}}{\prod
}(1-e^{-p\alpha }t^2)}.
$$
Here $\{e^\alpha \}$ is the standard notation
for the $X$-grading, and $t$ is the variable denoting the homological
degree.
$\qed$
\bs
The word ``minimal" is explained by the
following result of Ginzburg and Kumar.

\subsubsection{}
\Lemma
(see  [GK], Theorem 2.5.)
$$ ch(\ext_{{\u}^{-}\mod}\bul (k,k),t)=
\frac{\underset{w\in W}{\sum }e^{-w(\rho )+\rho }t^{l(w)}}
{\underset{\alpha \in R^{+}}{\prod }(1-e^{p\alpha }t^2)} \qed.
$$
\bs
One can construct  $R_{\min}^{-\bullet}(k)$
for $\u({\frak sl}_2)$ explicitly, and one easily obtains the
following statement.

\subsubsection{}
\Lemma
$$
ch(\underset{\lambda \in \Z}{\bigoplus }\Exts_{\C (\u ({\frak sl}_2))}
(k,k\langle p\lambda\rangle),t)=\frac{e^{p\alpha}(t+t^{-1})}
{(1-e^{p\alpha }t^2)(1-e^{p\alpha}t^{-2})}
$$
Here $\alpha$ is the only positive root of $\frak{sl}_2$.
$\qed$
\subsection{}
To obtain a character formula for  semiinfinite cohomology over other
finite  quantum groups we are going to use the filtrations $F'$ on
quantum groups defined in 2.4.6.  We  begin with the calculation of
$\underset{\lambda \in X}{\bigoplus } \Exts_{\C(\gr^{F'}(\u
))}(k,k(p\lambda ))$, using the following statement.

\subsubsection{}
\Lemma There exists a $\gr({\u}^{+})$-free resolution of
the trivial $\gr^{F'}(\u )$-module $k$ with a space of
$\gr({\u}^+)$-generators
$\underset{\alpha \in R^{+}}
{\bigotimes }\Lambda (\xi^\alpha )\otimes
\underset{\alpha \in R^{+}}
{\bigotimes }S(\eta^{p\alpha } ).$
Here $\{\xi ^\alpha \}$ denotes the exterior algebra generators of the
homological degree $-1$ and the $X$-grading $\alpha $, $\{\eta
^{p\alpha }\}$ denotes the symmetric algebra generators of  the
homological degree $-2$ and the $X$-grading $p\alpha$.
\ms
\dok
Consider the subalgebras $k(x_{\alpha})$, $\alpha\in R^+$, in
$\gr(\u^+)$ generated by the images of the root elements from $\u^+$.
Each algebra $k(x_{\alpha})$ is isomorphic to
$k[x_{\alpha}]/(x_{\alpha}^p)$, and $\gr^{F}(\u^+)$ is the twisted
tensor product of the algebras $k(x_{\alpha})$.  That is,
$\gr^F(\u^+)\cong \underset{\alpha \in R^+}{\bigotimes}k(x_{\alpha})$
as a vector space, and the following relations are satisfied:
$$ (*)
\hspace{3cm}
x_{\alpha}x_{\beta}=\zeta^{(\alpha|\beta)}x_{\beta}x_{\alpha}
$$
when
$\alpha<\beta$ in the chosen convex ordering on the set of positive
roots.

Moreover as proved by Lemma~\ref{filtration}  elements from $\u^-$
commute with elements of $\gr^F(\u^+$ in $\gr^{F'}(\u)$.  Thus it is
enough to construct $R_{\min}^{+\bullet}(k)$ for every algebra
$k(x_{\alpha})$ and to take tensor product of these resolutions over
the set of positive roots.

For the latter algebra the required
resolution looks as follows:
$$
0\map k(x_\alpha )x_\alpha^p\map k(x_\alpha )x_\alpha
\map k(x_{\alpha})\overset{\epsilon}{\map} k\map 0.
$$
The algebra
$gr^{F}({\u}^+)$  acts on the tensor product of such complexes by the set of
positive roots by the commutation rule (*),
the action of $\u^0$  on it comes from the $X$-grading and $\u^-$ acts
on the complex by zero.$\qed$ \bs We denote this resolution by
$R^{+\bullet}_{\min ,\gr^{F'}({\u})}(k)$.

\subsubsection{}
\Prop
$$
ch\left( \underset{\lambda \in X}
{\bigoplus }\Exts_{\C(\gr^{F'}(\u ))}(k,k\langle p\lambda \rangle),t\right)
=t^{-\dim{\frak n}^-}\frac {e^{2p\rho } \underset{w\in W}{\sum}
t^{2l(w)}}{\underset{\alpha \in R^+}{\prod }(1-e^{p\alpha
}t^2)(1-e^{p\alpha }t^{-2})}.
$$
Here as before $t$ denotes the
homological grading, $\{e^\alpha \}$ is the standard notation for the
$X$-grading. The equality is understood as an equality of power series in
variables numbered by the set of the generators of $X$ with coefficients
in  $k[t,t^{-1}]$.
\ms
\dok
$\gr^{F'}(\u^+)$ is a
Frobenius algebra, hence $\gr^{F'}(\u )^{\sharp}=\gr^{F'}(\u )$,
and
$$
\underset{\lambda \in X}{\bigoplus }\Exts_{\C(\gr^{F'}(\u ))}
(k,k\langle p\lambda \rangle) =H\bul (\Hom_{\gr^{F'}(\u )}\bul
(P^{-\bullet}(k),R_{\min ,\gr^{F'}(\u )}^{+\bullet}(k))).
$$
Here $P^{-\bullet}(k)$ is an
arbitary ${\u}^-$-free resolution  of
$k$ belonging to $\dC$.  The fact that
$\gr^{F'}({\u}^+)$ is Frobenius also implies that
$R_{\min ,\gr^{F'}( {\u}) }^{+\bullet}(k)$
consists of $\gr^{F'}(\u^+)$-cofree modules with the space of
cogenerators
$$
\underset{\alpha \in R^{+}}{\bigotimes }\Lambda (\xi ^\alpha )
\otimes \underset{\alpha \in R^+}{\bigotimes }
S(\eta ^{p\alpha })\otimes k_{( p-1)2\rho }.
$$
Consider the spectral sequence of the bicomplex
$$
\Hom_{\gr^{F'}(\u )}\bul (P^{-\bullet}(k),R_{\min
,\gr^{F'}(\u )}^{+\bullet}(k)).
$$
The term $E_1$ looks as follows:
$$
\underset{\lambda\in X}{\bigoplus}\ext_{\C(\b^-)}\bul(k,\underset{\alpha \in
R^+}{\bigotimes } \Lambda (\xi ^\alpha )\otimes k_{(p-1) 2\rho}\langle
p\lambda\rangle) \otimes \underset{\alpha \in R^+}{\bigotimes} S(\eta
^{p\alpha }),
$$
where the space $\underset{\alpha \in
R^{+}}{\bigotimes }\Lambda (\xi ^\alpha )\otimes k_{(p-1)2\rho }$ is
the direct sum of one dimensional $\b^-$-modules.

\subsubsection{}
\Lemma (see [GK], Theorem
2.5.)\\                                     \label{GK}

\qquad(i) $\underset{\lambda \in X}{\bigoplus }\ext_{\C(b^-)}\bul(k,
k_\mu\langle p\lambda \rangle)=0$ when $\mu \neq w\left( \rho \right) -\rho$;

\qquad(ii) for $\mu =w( \rho )-\rho$
$$
ch\left( \underset{\lambda \in X}{\bigoplus }\ext_{\C(\b^-)}\bul(k,
k_{\mu} \langle p\lambda \rangle),t\right) =
\frac {t^{l( w)}} {\underset{\alpha\in R^+}{\prod}(1-e^{p\alpha }t^2) }.
\qed
$$
\subsubsection{}
\Lemma                                                \label{jantzen}
(see   [J], part II, Lemma 12.10) Let $\alpha _1,\ldots ,\alpha
_k$ and $\beta _1,\ldots ,\beta _k$ be the two sets of pairwise distinct
positive roots such that
$$
\alpha _1+\ldots +\alpha
_k\equiv \beta _1+\ldots +\beta _mmod\left( pX\right).
$$
Then for $p>2(h-1)$
$$
\alpha _1+\ldots +\alpha _k=\beta _1+\ldots
+\beta _m.
$$
Here $h$ denotes the Coxeter number of the  root
system (see  [J], p.262).$\qed$
\bs
For any element of the Weyl group
$w\in W$ there exists a unique element $w'\in W$, such that
$w(\rho )+w'(\rho )=0$. Then $w(\rho )-\rho =-2\rho - (w'(\rho )-\rho )$.
In [GK], 2.5 it is proved that for every $w'\in W$
$$
\dim\ _{\rho-w'(\rho)}\left(\underset{\alpha\in R^+}
{\bigotimes}\Lambda(\xi^\alpha)\right)=1.
$$
Thus for every $w\in W$
$$
\dim\ _{w(\rho)-\rho}\left(\underset{\alpha\in R^+}
{\bigotimes}\Lambda(\xi^\alpha)\otimes k_{-2\rho}\right)=1.
$$
Noting that $l(w')=\dim{\frak n}^--l(w)$ if $w(\rho)+w'(\rho)=0$ we
see that all the nonzero entries of the term $E_1$ of the spectral sequence
are of the same
parity whence it degenerates, and the proposition follows from Lemmas~\ref{GK}
and~\ref{jantzen}.
$\qed$

\subsection{}
\label{answer}
\Theorem
$$
ch\left( \underset{\lambda\in X}{\bigoplus }
\Exts_{\C}(k,k\langle p\lambda \rangle),t\right)
=\frac{t^{-\dim{\frak n}^-}e^{2p\rho }\underset{w\in W}{\sum }
t^{2l(w)}}
{\underset{\alpha \in R^+}{\prod }(1-e^{p\alpha }t^2)
(1-e^{p\alpha }t^{-2})}.
$$
\ms
\dok
The spectral sequence arising from  the filtration of
the complex $\stand(k,k)$, which is induced by the filtration $F'$ on
${\u}$, degenerates as all the nonzero semiinfinite cohomology
of the trivial module over $\gr^{F'}(\u )$
are of  the same parity.$\qed$

\section{$\protect\Tors^{\C}(L,M)$ and conformal blocks}

In this section we are going to compare our semiinfinite cohomology
of quantum groups with the functor defined by Finkelberg and
Schechtman.

\subsection{}
First we  define   the semiinfinite homology of quantum groups.

\subsubsection{}
\Def For $M\in {\C}$, $L\in {\C}^{(r)}$  we set
$$
\Tors^{\C}(L,M)=\ _0(H\bul (P_{\u^-}\bul (L)\otimes_{\u}P_{\u^+}\bul(M))).
$$
Here $P_{\u^-}\bul (L)\in \C^{(r)\downarrow}$ is a
${\u}^-$-free resolution of $L$,
$P_{\u^+}\bul(M)\in \dC$ is a
$\u^+$-free resolution of $X$; and $_0(H\bul(\ldots))$
denotes the zeroth $X$-grading component.

Like in the third section, one can easily check that
the definition doesn't depend on the choice of
resolutions. Alternatively, this is
a corollary of the following comparison Lemma:

\subsubsection{}
\Lemma \label{compare} $\Exts_{\C}(L,M)=
\tor_{\frac \infty 2-\bullet}^{\C}(M^*,L)^*.$
\ms
\dok
The statement  follows immediately from the
definition of semiinfinite homology and the standard relation between
$\Hom_{\u}$ and $\otimes_{\u}$.
$\qed$
\bs
In particular semiinfinite $\tor$ is well defined as a functor on the
corresponding derived categories of $\u$-modules (see \ref{derived}).
We will need the following reformulation of
the statement of Lemma~\ref{main}(iii)
on the language of semiinfinite $\tor$ functor.

\subsubsection{} \label{main2}
\Lemma
For $L\bul\in\C^{(r)\downarrow}$, $M\bul\in\dC$
if $M\bul$ is both $\u^+$-projective
and $\u$-injective relative to $\u^-$ then
$\Tors^{\C}(L\bul,M\bul)=\ _{(0)}(H\bul(L\bul\otimes_{\u}M\bul)).$
$\qed$
\bs
Note that since both algebras $\u$ and $\u^+$ are Frobenius (see~\ref{frob}),
the condition of the previous lemma means simply that $M\bul\in\dC$
consists of $\u$-projective modules.

\subsubsection{}
\Rem
M.~Finkelberg and V.~Schechtman gave a geometric definition of semiinfinite
$\tor$ functor in the category $\C$ (see [FiS], part IV).
In [Ar] and [FiS], part IV, it is proved that the geometric definition
coincides with the one presented here.

\subsection{}
Semiinfinite homology plays an important part in the calculation of conformal
blocks. Recall  the definition of conformal blocks (see e.g. [A]).

Suppose for simplicity that our Cartan matrix is symmetric.
 Let  $\gamma \in R^{\vee}$ be the highest
root, and denote by $\Delta \in X$ the first  alcove,
$$
\Delta=\{ \lambda \in X|\langle \gamma,\gamma +\rho\rangle <p,\langle
\alpha_i^{\vee},\lambda+\rho\rangle >0,i=1,\ldots ,r\}.
$$

For $\lambda_1,\ldots ,\lambda_n \in \Delta$ the
conformal blocks  $\langle
L_{\lambda_1},\ldots ,L_{\lambda_n}\rangle$ are defined as the
maximal trivial direct summand
in $L_{\lambda_1}\otimes \ldots \otimes
L_{\lambda_n}$.

\subsection{}
\Theorem \label{embed} There exists a
natural inclusion
$$
\varphi :\ \langle L_{\lambda_1},\ldots ,L_{\lambda_n}\rangle \hookrightarrow
\tor_{\frac \infty 2+0}^{\C}(k,L_{\lambda_1}\otimes \ldots \otimes
L_{\lambda_n}\otimes L_{(p-1)2\rho}).
$$
The proof  follows immediately from the definition of the conformal blocks
and the following statement.

\subsubsection{}
\Prop
$\tor_{\frac \infty 2+0}^{\C}(k,L_{(p-1)2\rho})=k.$
\ms
\dok
Using~\ref{compare}
we are going to prove the corresponding
statement for semiinfinite cohomology.

Choose left resolutions of
$k$ and $L_{(p-1)2\rho}$ beginning  with
$M_0^+$ and $M_{-(p-1)2\rho}^-$ respectively and
satisfying  the conditions (U) and (D) respectively.
Such resolutions exist, for example one can take the standard resolution
$\Barb(\u,\b^-,\overline{M}_0^+)$ of the kernel of the canonical projection
$M^+_0\map  k$
(resp.  the standard resolution $\Barb(\u,\b^+,\overline{M}_{(p-1)2\rho}^-)$
of the kernel of the canonical projection
$M_{(p-1)2\rho}^-\map
L_{(p-1)2\rho}$).  Denote these resolutions  by
 $P^{(-)\bullet}(L_{(p-1)2\rho)}$ and
$P^{(+)\bullet}(k)$ respectively.

Using~\ref{verma}(iii) we see that
$\Hom_{\C}(P^{(-)n}(L_{(p-1)2\rho}),P^{(+)m}(k))$ is nonzero only
when $n=m=0$, and
$$
\Hom_{\C}(P^{(-)0}(L_{(p-1)2\rho}),P^{(+)0}(k))=k.\qed
$$

\subsubsection{}
\label{example} The inclusion $\varphi$ in general is not bijective.
Consider the following example:  ${\frak g}={\frak sl}_2,p=5,n=4$,
$\lambda_1=\lambda_2=2,$ $\lambda_3=\lambda_4=3$ (we have identified $X$
with $\Z$, so $\rho=1$, and $(p-1)2\rho=8$).

\subsubsection{}
\Lemma $\tor_{\frac \infty 2+0}^{{\C}}(k,P_0)=k$, where $P_0$ is the
projective covering of the trivial module.
\ms
\dok
The statement  follows easily from~\ref{main2}.$\qed$
\bs
As the maximal trivial direct summand of $P_0$ is zero, it is enough to find
this module among the direct summands of $L_2\otimes L_2\otimes L_3\otimes
L_3\otimes L_8$.  If we find a projective direct summand  $P$ in
$V=L_2\otimes L_2\otimes L_3\otimes L_3$ with the highest weight $0$, then
$V\otimes L_8$ will contain a projective direct summand with the
highest weight $8$, i. e.  $P_0$.

It is well known (see  [L2], Proposition 7.1), that all the modules
$L_{\lambda},\lambda \in \Delta,$ lift to the simple $U_3$-modules.  It
follows from the results of [A] that for ${\frak g}={\frak sl}_2$,
$\lambda_1,\ldots ,\lambda_n \in \Delta$, the module $L_{\lambda_1}\otimes
\ldots \otimes L_{\lambda_n}$ is a direct sum of projective $U_3$-modules
and simple $U_3$-modules with highest weights  in $\Delta$.

Thus the $U_3$-module $L_2\otimes L_2\otimes L_3\otimes L_3$ contains
the idencomposable projective $U_3$-module $\widetilde{P}_{8}$
(which is the projective covering of the simple $U_3$-module
${\widetilde{L}}_{8}$) with the highest weight $10$ as an
$U_3$-direct summand.  One can check easily that when restricted to
${\u}$ the module $\widetilde{P}_{8}$ contains as a direct
summand $P_{-2}$ --- the projective covering of the ${\hbox{\bf
u}}$-module $L_{-2}$. But the highest weight of $P_{-2}$ is $0$.

We conclude that $L_2\otimes L_2\otimes L_3\otimes L_3\otimes L_8$
contains a direct summand $P_0$. Hence the semiinfinite $\tor$ is strictly
bigger than  conformal blocks.

\subsection{}
Now we construct a certain duality on semiinfinite $\tor$ spaces which
corresponds to Poincar\`e duality in Finkelberg-Schechtman
interpretation.

\subsubsection{}
Denote by $\widetilde{\u}$ the finite
quantum group defined in the same way as
$\u$, but with $\zeta$ replaced by $\zeta^{-1}$ in the defining
relations.

\subsubsection{}
\Lemma The map
$
\phi:\ E_i\mapsto F_i,\ F_i\mapsto E_i,\ K_i\mapsto K_i^{-1}
$
defines an antiisomorphism of algebras
$\phi:\text{ }{\u}\map  {\widetilde{\u}}$.
$\qed$
\bs
Denote by ${\widetilde{\C}}$ the category of
left ${\widetilde{\u}}$-modules satisfying the conditions of the
type~\ref{wedge}.  Define the functor
$$
D:\  {\C}\map {\widetilde{\C}}^{opp},\ D(M)=\Hom_k(M,k)
$$
with the natural left action of ${\widetilde{\u}}$ constructed as follows:
$$
u\cdot f(m):=f(\phi^{-1}(u)m),\text{ where }u\in {\widetilde{\hbox{\bf
u}}},f\in \Hom_k(M,k).
$$
One can check directly that the supports of the
modules $M\in {\C}$ and $D(M)\in {\widetilde{\C}}$
coincide.

\subsubsection{}
\Prop
There exists a nondegenerate pairing
$$
\langle \ ,\ \rangle :\ \Tors^{\C}(k,M)\times \tor_{\frac \infty
2-\bullet}^{\widetilde{\C}}(k,D(M))\map  k.
$$
\dok
By~\ref{compare}
$$
\Tors^{\C}(k,M)=(\ext_{\C}^{\frac\infty2-\bullet}(M,k))^*,\
\Tors^{\widetilde{\C}}(k,D(M))=(\ext_{\widetilde{\C}}^{\frac\infty2-\bullet}
(D(M),k))^*,
$$
thus it is enough to construct a nondegenerate pairing
$$
\Exts_{\C}(M,k)\times\ext_{\widetilde{\C}}^{\frac\infty2-\bullet}(D(M),k)\map k.
$$
Choose a resolution $R\bul(k)\in\dC$ consisting of projective $\u$-modules.
Then by~\ref{main2} $\Exts(M,k)=\ _0(H\bul(\Hom\bul_{\u}(M,R\bul(k))))$,
and by definition of $D$
$$
\Exts_{\widetilde{\C}}(D(M),k)=\ _0\left(H\bul\left(
\Hom\bul_{\widetilde{\u}}(D(M),D(R\bul(k))\right)\right)=\
_0\left(\Hom\bul_{\u}(R\bul(k),M)\right).
$$
We may assume that $R\bul(k)$
consists of modules of the form $L=\Coind_{\u^0}^{\u}(V)$.
For such a module consider the canonical isomorphism
\begin{multline*}
  \mu:\ _0(\Hom_{\u}(M,L))\widetilde{\map}\ _0\left(\Hom\bul_{\u^0}(M,V)\right)
  \widetilde{\map}\ _0\left(\Hom\bul_{\u^0}(V,M)\right)^*\widetilde{\map}\\
  _0\left(\Hom\bul_{\u}(\Ind_{\u^0}^{\u}(V),M)\right)^*\widetilde{\map}
  \ _0\left(\Hom\bul_{\u}(L,M)\right)^*
\end{multline*}
The last equality uses the fact
that $\Ind_{\u^0}^{\u}(V)=\Coind_{\u^0}^{\u}(V)$ since
all the algebras $\u^+$, $\u^-$ and $\u$ are Frobenius.

We are to check that $\mu$ commutes with  morphisms of $\u$-modules
$L_1\map L_2$.  But the general statement follows easily from the
statement for $L_1$ and $L_2$ being $\u$-free. So we are to check that
$\mu:\ \Hom_{\u}(\u,M)^*\widetilde{\map} \Hom_{\u}(M,\u)$ is an
isomorphism of right $\u$-modules. This can be verified directly.  So
$\mu$ induces a nondegenerate pairing of the complexes:
$$
\langle\
,\ \rangle:\
\Hom_{\u}\bul(R\bul(k),M)\times\Hom_{\u}\bul(M,R\bul(k))\map k.
$$
that becomes the required pairing on the semiinfinite $\tor$ spaces.
$\qed$

\appendix
\section{Distributions on the nilpotent cone}

In this section we are going to give a geometric interpretation of the
 character of the semiinfinite cohomology of the trivial
module over the finite quantum group.

\subsection{}
Consider the nilpotent cone
${\N}\subset {\g}$.
${\N}$ is a singular affine algebraic variety containing the positive
nilpotent subalgebra ${\frak n}\subset {\N}\subset {\g}$.  Denote
by ${\cal F}({\N})$ the space of complex algebraic functions on ${\cal
N}$.  The adjoint action of $\frak g$ preserves $\N$ and induces a
representation of $\frak g$ in  ${\cal F}({\N})$.
Our considerations are parallel to the following result due to
Ginzburg and Kumar [GK].

\subsubsection{}
{\bf Proposition:} (see [GK], Theorem 5)
(i) $\underset{\lambda \in
X}{\bigoplus }\ext_{\C}^{2n+1}( k,k ( p\lambda )
) =0$ for any $n\ge 0$;\\
(ii) $\underset{\lambda \in X}{\bigoplus
}\ext_{\C} ^{2n}(k,k ( p\lambda )) = {\cal
F}^{n}({\N})$ as ${\frak g}$-modules. Here the grading on the right
hand side is the grading by homogeneous degree of functions on ${\cal
N}.\qed$
\bs
B.~Feigin has proposed the following conjecture.

\subsubsection{}
\Con
The $\frak g$-module
$\underset {\lambda \in X}{\bigoplus} \Exts_{\C}(k,k\langle p\lambda\rangle)$
is isomorphic to the $\frak g$-module of distributions on $\N$
with support in ${\frak n}\subset {\N}$.$\qed$
\bs
To formulate the exact statement we will need
several well known facts about  the geometry of the nilpotent cone.

We will need the Grothendieck-Springer resolution of ${\N}$.  Choose a
maximal torus $T\subset G$ and  a Borel subgroup $T\subset B\subset
G$. Consider the flag variety ${\B}=G/B$. Le ${\frak b}\subset {\frak
g}$ be the Lie algebra of $B$, and let $\frak n$ be its nilpotent
radical.

\subsubsection{}
\Lemma (see [CG], 3.1.36)
(i) The natural map
$$
\sigma :T^*({\B})\map  {\N},
$$
where $T^*({\B})$ is the cotangent bundle to ${\B}$, is a resolution of
singularities of $\N$.\\
(ii)
$$
\sigma ^{-1}({\frak n}) = \underset {w\in W}{\coprod
}T^*_{\C_w}{\B}\subset T^*(\B ),
$$
where $T^*_{\C_w}{\B}$ denotes the conormal bundle to a $B$-orbit ${\cal
C}_w$.
$\qed$
\bs
Denote the union of the conormal bundles to the $B$-orbits ${\C}_w$
by ${\cal S}$.

\subsection{}
We use some results and methods
due to Kempf [K].  The notion of cohomology with support in a locally
closed subvariety was  investigated in that paper.  Kempf
also introduced the action of a Lie algebra on local cohomology in the
case when the corresponding Lie group acts on the ambient space.

We  denote the local cohomology  of $X$ with support in $Y\subset X$ and with
coefficients in a sheaf ${\cal F}$ by $H_Y\bul (X,{\cal
F})$.

\subsubsection{}
\Lemma

\quad(i) $\frak g$ acts naturally on
$H_{\frak n}\bul ({\N},{\O}_{\N})$. Here ${\O}_{\cal
N}$ denotes the structure sheaf of $\N$.

\quad(ii) There is a natural isomorphism
of $\frak g$-modules
$H_{\frak n}\bul ({\N},{\O}_{\N})\widetilde{\map }
H_{\cal S}\bul (T^{*}({\B}),{\cal
O}_{T^{*}(\B )})$.
\ms
\dok
(i) Follows from [K], Lemma 11.1.

(ii) The existence of the isomorphism
of vector spaces in question follows from the fact that $R^0\sigma_*
{\O}_{T^*{\B}}={\O}_{\N}$, and $R^i\sigma_*
{\O}_{T^*{\B}}=0$. The isomorphism commutes with
the $\frak g$-action since $\sigma$ is $G$-equivariant.
$\qed$
\bs
We will calculate the character of the $\frak g$-module
$H\bul_{\cal S}(T^*({\B}),{\O}_{T^*({\B})})$.
Denote by $s$ the homological grading. The
grading by the homogeneous degree is denoted by $t$.

\subsubsection{}
\Theorem
$$
ch\left(H_{\cal S}\bul (T^*({\B}),\O_{T^*({\B})}),t\right)=
e^{2\rho}\underset {w\in W}{\sum }\frac
{t^{l(w)}}{\underset {\alpha \in
R^+}{\prod}(1-e^{\alpha}t)(1-e^{\alpha}t^{-1})}s^{\dim{\frak n}^-},
$$
the equality here is the equality of power series in the variables
numbered  by the generators of $X$ with coefficients in
$k[t,t^{-1}]$.
\ms
\dok
The statement follows immediately
from the  Lemmas~\ref{a},\ref{b} and~\ref{c}.
$\qed$
\bs
Comparing the answer with  Theorem~\ref{answer} we obtain the following fact.

\subsubsection{}
\Cor
Up to a shift of grading  in $t$ the character of the
semiinfinite cohomology of the trivial module over the finite quantum
group coincides with the character of $H^{\bullet}_{\frak n}({\N},{\O}_{\N})$.
$\qed$

\subsubsection{}
\Lemma                                      \label{a}
$ch\left(H_{\cal S}^{\bullet}(T^*({\cal
B}),{\O}_{T^{*}({\B})}),t\right) = \underset {w\in W}{\sum }
ch\left(H_{T^*_{{\C}_{w}}{\B}}^{\bullet}(T^*({\B}),{\cal
O}_{T^*({\B})}),t\right).$
\ms
\dok
Fix a total linear order on the set $\{{\cal C}_w|w\in W\}$ compatible
with the natural partial order by inclusion.

We introduce a the filtration on $\cal S$ by subspaces ${\cal
S}_w:=\underset {l(w')<w}{\bigcup}T^*_{{\C}_{w'}}{\B}$.

We prove by induction that
$$
ch\left(H\bul_{{\cal S}_w}
\left(T^*({\cal B}),{\cal O}_{T^*({\cal B})}\right),t\right)=
\underset{w'\le w}{\sum}
ch\left(H\bul_{T^*({\cal C}_{w'})}
\left(T^*({\cal B}),{\cal O}_{T^*({\cal B})}\right),t\right).
$$
For $w=1,\ {\cal C}_w=pt,$ the statement is evident.

If $w$ follows $w'$ directly in the total linear order then
${\cal S}_{w'}$ is closed in
${\cal S}_w$ and ${\cal S}_{w}\backslash{\cal S}_{w'}
=T^*({\cal C}_w)$. Then by definition of
local cohomology there exists a long exact sequence
\begin{multline*}
\ldots\map H^i_{{\cal S}_{w'}}
\left(T^*({\cal B}),{\cal O}_{T^*({\cal B})}\right)
\map H^i_{{\cal S}_w}\left(T^*({\cal B}),{\cal O}_{T^*({\cal B})}\right)
\map\\ \map H^i_{{\cal C}_w}\left(T^*({\cal B}),{\cal O}_{T^*({\cal B})}\right)
\map H^{i+1}_{{\cal S}_{w'}}\left(T^*({\cal B}),{\cal O}_{T^*({\cal B})}\right)
\map\ldots
\end{multline*}
The cohomology of a nonsingular variety with support in a nonsingular (locally
closed) subvariety is nonzero only in one degree equal to the codimension
of the subvariety, thus the statement is proved.
$\qed$

\subsubsection{}
\Lemma \label{b}(see [K], 6.4) There exists a $T$-equivariant
neighbourhood of a $B$-orbit $P({\C}_{w})$ such that\\
(i)
$$
P({\cal
C}_{w})\cong \underset {\alpha \in R^+,w(\alpha) \in R^-}{\prod}V(\rho
-w(\rho ))\times \underset {\alpha \in R^-,w(\alpha )\in R^-}{\prod
}V(\rho -w(\rho ))
$$
as topological spaces with the action of $T$. Here
$V(\lambda)$ denotes the one dimensional representation of $T$ of   weight
$\lambda$;\\
(ii) under this identification ${\C}_{w}$ corresponds to
the first factor.$\qed$

\subsubsection{}
\Lemma          \label{c}
$$
\underset
{w\in W}{\sum }t^{l(w)}\left( \underset {\alpha \in R^+,w(\alpha)\in
R^+}{\prod }(1-e^{\alpha}t)\underset {\alpha \in R^+,w(\alpha)\in
R^-}{\prod}(1-e^{\alpha}t^{-1})\right)
$$
$$=\left( \underset {\alpha \in
R^+}{\prod}(1-e^{\alpha}t)\right)\left(\underset{w\in
W}{\sum}t^{l(w)}\right).\qed
$$

\section{Semiinfinite cohomology as a derived functor}

Here we discuss the notion of a K-semijective complex introduced by Voronov
(see [V]) and relations between our approach to semiinfinite cohomology of
associative algebras and Voronov's investigation of Lie algebras' semiinfinite
cohomology.

It turns out that Voronov's resuts remain true not only in the case of
graded Lie elgebras but also in a more general setting af an
associative algebra $A$ equipped with a triangular decomposition
$A=B\otimes N$. So we remain in the situation of~\ref{setup}.

We begin with recalling several basic results of semiinfinite
homological algebra.

\subsection{}
Let ${\cal O}(A)$ be the category of $X$-graded $A$-modules
$M=\underset{\lambda\in X}{\bigoplus}\ _{\lambda}M$
satisfying the following condition:
there exist $\beta_1,\ldots,\beta_s\in X$ such that $\supp M\subset
\bigcup_{i=1}^sX^+_\Q(\beta_i)$, with morphisms being morphisms of $A$-modules
that preserve $X$-gradings. The corresponding category of $X$-graded
$N$-modules is denoted by ${\cal O}(N)$.

We denote the homotopy category of unbonded complexes over ${\cal
O}(A)$ (resp., over ${\cal O}(N)$) with morphisms being morphisms of
complexes modulo nillhomotopies, by ${\cal K}(A)$ (resp., by ${\cal
K}(N)$).

\subsubsection{}
\Def (see [V], 3.3) An object $S\bul\in {\cal K}(A)$ is called
K-semijective if

(i) it is K-projective over $N$, i.~e.
$
\Hom_{{\cal K}(N)}(S\bul,V\bul)=0
$
for any (possibly unbounded) complex $V\bul\in {\cal K}(N)$;

(ii) it is K-injective over $A$ relative to $N$, i.~e.
$
\Hom_{{\cal K}(A)}(V\bul,S\bul)=0
$
for any complex $V\bul\in {\cal K}(A)$ such
that $\Hom_{{\cal K}(N)}(V\bul,V\bul)=0$
(in partiqular $V\bul$ is acyclic).

Note that
our definition is ``turned upside down'' with respect to the Voronov's one.

Note also that $A$-modules that are both $N$-projective and
$A$-injective relative to $N$ (see~\ref{rel}) evidently are
K-semijective. We call such $A$-modules semijective (without K).

\subsubsection{} \label{bound}
\Lemma (i) Any bounded from above complex of $N$-projective $A$-modules
is K-projective over $N$;

(ii) any bounded from below complex of $A$-injective relative to $N$ modules
is K-injective relative to $N$;

(iii) in particular any bounded complex of semijective modules is K-semijective.
$\qed$

The following statement is the main achievement of semiinfinite
homological algebra.

\subsection{}
\Theorem \label{sem}
(see [V], Theorem 3.3)
Let ${\cal K}({\cal SJ}(A))$ be the homotopical
catregory of K-semijective
complexes over ${\cal O}(A)$. ${\cal D}(A)$ denotes the unbounded
derived category of complexes over ${\cal O}(A)$. Then the functor of
localization by the class of quasiisomorphisms provides a natural
equivalence of triangulated categories
$$
{\cal K}({\cal
SJ}(A))\widetilde{\map}{\cal D}(A).\qed
$$

\subsubsection{}
Consider the following resolution $R\bul(M)$ of an $A$-module $M\in
{\cal O}(A)$:
$$
R\bul(M):=\Hom_A\bul(\Barb(A,N,A),\Barb(A,B,M)),
$$
where $A$ in $\Barb(A,N,A)$ is considered as a $A-A$~bimodule.
Evidently $R\bul(M)\in \upC(A)$.

\subsubsection{}
\Lemma
 \label{main3}
$R\bul(M)\in\upC(A)$ is K-semijective.
\ms
\dok
First note that $\Barb(A,N,A)$ is homotopically equivalent to $A$ as a
$A-N$~bimodule --- the homotopy is provided by the map
$$
a_0\otimes\ldots\otimes a_n\otimes a\map a_0\otimes\ldots\otimes a_n
\otimes a\otimes 1.
$$
Thus $R\bul(M)$ is homotopically equivalent to
$$
\Hom_A\bul(A,\Barb(A,B,M))=\Barb(A,B,M)
$$
as a $N$-module. In partiqular for any exact complex $V\bul\in{\cal O}(A)$
$$
\Hom_{{\cal K}(N)}(R\bul(M),V\bul)=\Hom_{{\cal K}(N)}(\Barb(A,B,M),V\bul)=0.
$$
Thus $R\bul(M)$ is K-projective over $N$.
To prove that $R\bul(M)$ is $A$-injective relative to $N$ note that
for a complex of $A$-modules $V\bul$ homotopically equivalent to zero over
$N$
\begin{multline*}
\Hom\bul_A(V\bul,\Hom\bul_A(\Barb(A,N,A),\Barb(A,B,M)))=\\
\Hom\bul_A(\Barb(A,N,A)\otimes_AV\bul,\Barb(A,B,M))=\\
\Hom\bul_A(\Barb(A,N,V\bul),\Barb(A,B,M)).
\end{multline*}
Since $V\bul\cong 0$ in ${\cal K}(N)$,  each line $\Bar^p(A,N,V\bul)$
of the bicomplex $\Barb(A,N,V\bul)$ is homotopically equivalent  to zero
over $A$. As $\Bar^p(A,N,V\bul)\ne 0$ only for $p\le 0$, the total complex
of $\Barb(A,N,V\bul)$ is also homotopically equivalent to zero over $A$.

Thus $\Hom_{{\cal K}(A)}(V\bul,\Hom_A\bul(\Barb(A,N,A),\Barb(A,B,M)))=0$,
and we are done.
$\qed$

Thus by \ref{main}(iii) one
can treat $\Exts_A(L,M)$ as an exotic derived functor
of the functor $M\mapsto \Hom_{\oppA}(L,S_A\otimes_AM)$ (cf. [V], 3.9).

In particular in [V], 3.2.1 it is proved that in the case of $A=U({\frak a})$
for some graded Lie algebra $\frak a$ the algebra $\oppA$ differs from $A$
by a 2-cocycle of the Lie algebra $\frak a$. One can check directly
that the functor $\Hom_{\oppA}(k,S_A\otimes_A*)$ coincides with the functor of
semiinvariants defined in [V], 3.6.

\subsubsection{}
Recall the construction of the standard complex for the computation of
Lie algebra semiinfinite cohomology.

For a graded module $M$ over a graded Lie algebra
$\frak{a}=\underset{n\in\Z}{\bigoplus}{\frak a}_n$
the standard resolution with respect to the graded Lie subalgebra
${\frak b}\subset{\frak a}$  looks as follows:
$$
\St({\frak a},{\frak b},M):=\left(U({\frak a})\otimes_{U({\frak b})}
\Lambda\bul({\frak a}/{\frak b})\right)\otimes M.
$$
Here the $\frak b$-module $\Lambda({\frak a}/{\frak b})$ is just the direct
sum of the exterior powers of the ${\frak b}$-representation in
${\frak a}/{\frak b}$,  tensor product of $\frak a$-modules over the base field
is defined using the Hopf algebra structure on $U({\frak a})$,
the differential is written as follows:
\begin{multline*}
\operatorname{d}((u\otimes
\overline{a}_1\wedge\ldots\wedge\overline{a}_n)\otimes m
=\sum_{i=1}^n(-1)^i
(ua_i\otimes\overline{a}_1\wedge\ldots\wedge\overline{a}_{i-1}
\wedge\overline{a}_{i+1}\wedge\ldots\wedge\overline{a}_n)\otimes m\\
+\sum_{i<j}(-1)^{i+j}(u\otimes\overline{[a_i,a_j]}\wedge\overline{a}_1\wedge
\ldots\wedge\overline{a}_{i-1}\wedge\overline{a}_{i+1}\wedge\ldots\wedge
\overline{a}_{j-1}\wedge\overline{a}_{i+1}
\wedge\ldots\wedge\overline{a}_n)\otimes m.
\end{multline*}

Here $\overline{a}_i\in{\frak a}/{\frak b}$. One can check that the
differential is correctly defined.

Consider the triangular decomposition of the Lie algebra $\frak a$:
$$
{\frak a}^{\le 0}:=\underset{n\le 0}{\bigoplus}{\frak a}_n,\
{\frak a}^{>0}:=\underset{n>0}{\bigoplus}{\frak a}_n,\
{\frak a}={\frak a}^{\le 0}\oplus{\frak a}^{>0}
$$
as a vector space.

Clearly $\St({\frak a},{\frak a}^{\le 0},M)$ belongs to $\upC(U({\frak a}))$.

Let ${\frak a}^{\sharp}$ be the central extension of ${\frak a}$ such that
$U({\frak a})^{\sharp}=U({\frak a}^{\sharp})$.
Then, as before,
$$
K\bul_{U({\frak a})}(k,M)
=\Hom\bul_{U({\frak a}^{\sharp})}
(\St({\frak a}^{\sharp},{\frak a}^{\sharp>0},k),
S_{U({\frak a})}\otimes_{U({\frak a})}\St({\frak a},{\frak a}^{\le 0},M)),
$$
$$
\Exts_{U({\frak a})}(k,M)=H\bul(K\bul_{U{\frak a})}(k,M)).
$$
The  complex $K\bul_{U({\frak a})}(k,M)$ is exactly the standard complex for
the computation of Lie algebra semiinfinite cohomology
consisting of semiinfinite exterior powers (see e.g. [V], 2.5).
That gives another proof of coincidence
of semiinfinite $\ext$ functor and Lie algebra semiinfinite cohomology.

\section*{References}

[A] H.H.Andersen. {\it Representations of quantum groups, invariants of
3-manifolds and semisimple tensor categories.} Israel Math. Conf. Proc.
v.\hbox{\bf 7}, (1993), p.1-12.\\
$\text{[AJS] H.H.Andersen}$, J.C.Jantzen, W.Soergel. {\it Representations
of quantum groups at p-th roots of unity and of semisimple groups in
characteristic p: independence of p.} Asterisque \hbox{\bf 220}, (1994).\\
$\text{[Ar]}$ S.M.Arkhipov. {\it Semiinfinite cohomology of quantum
groups at roots of unity.} Preprint, (1994).\\
$\text{[CG]}$ N.Chriss, V.Ginzburg. {\it Representation theory and complex
geometry.}  Boston, Birkh\"auser,  (1995).\\
$\text{[DCK]}$ C.De Concini, V.Kac, {\it Representations of quantum
groups at roots of unity}, in A.Connes et all  (eds.), Operator algebras,
unitary representations, enveloping algebras and invariant theory
 (Colloque Diximier),
 Proc. Paris 1989, (Progr. in math. \hbox{92}), Boston etc.,
Birkh\"auser, p.
471-506.\\ $\text{[DCKP] C.De Concini, V.Kac, C.Procesi.}$ {\it Some
remarkable degenerations of quantum groups.} Comm. Math. Phys.
\hbox{\bf 157},  (1993), p.405-427.\\ $\text{[F]}$ B.Feigin.
{\it Semi-infinite cohomology of Kac-Moody and Virasoro Lie algebras.} Usp.
Mat.  Nauk \hbox{\bf 33},no.2, 195-196, (1984),  (in Russian).\\
$\text{[FSV]}$ B.Feigin, V.Schechtman, A.Varchenko.  {\it On algebraic
equations satisfied by hypergeometric corellators in WZW models.}  II,
Comm. Math.  Phys. {\bf 170}, 219-247, (1995).     \\
$\text{[Fi]}$ M.Finkelberg. {\it An
equivalence of fusion categories.}  Harvard Ph.D. thesis, (1993).\\
$\text{[FiS]}$ M.Finkelberg, V.Schechtman.  {\it Localization of}
${\u}$-{\it modules. I. Intersection cohomology of real arrangements,}
 Preprint hep-th/9411050  (1994),1-23; II. {\it
 Configuration spaces and quantum groups,}
 Preprint q-alg/9412017  (1994),1-59; III  {\it
Tensor categories arising from configuration spaces.}
Preprint q-alg/9503013 (1995) 1-59; IV
{\it Localization on} $\mathbf{P}^1$, Preprint q-alg/9506011 (1995), 1-31.\\
$\text{[GeM]}$
S.I.Gelfand, Yu.I.Manin.  {\it Methods of homological algebra.} M:  Nauka,
 (1988).  (In Russian)\\ $\text{[GK]}$ V.Ginzburg, N.Kumar.  {\it Cohomology
of quantum groups at roots of unity.}  Duke Math. J. \hbox{\bf 69},  (1993),
179-198.\\ $\text{[H]}$ W.H.Hesselink. {\it On the character of the
nullcone.}  Math. Ann.  \hbox{\bf 252}, (1980), 179-182.\\ $\text{[J]}$
J.C.Jantzen.  {\it Representations of algebraic groups.}  Pure and Appl.
Math.  \hbox{\bf 131}, Orlando, Fla., Academic Press, (1987).\\
$\text{[K]}$ G.Kempf.  {\it The Grothendieck-Cousin complex of an induced
representation.}  Adv. in Math., \hbox{\bf 29}, 310-396,  (1978).\\
$\text{[KL 1,2,3,4]}$  D.Kazhdan and G.Lusztig. {\it Tensor structures
arising from affine Lie algebras.} I, J.  Amer.  Math. Soc. {\bf 6},
 (1993), 905-947; II, J. Amer. Math. Soc. {\bf 6}  (1993), 949-1011; III, J.
Amer. Math. Soc.  {\bf 7}  (1994), 335-381; IV, J. Amer. Math. Soc. {\bf 7},
 (1994), 383-453.\\ $\text{[L1]}$ G.Lusztig.  {\it Quantum groups at roots
of unity.}  Geom.  Dedicata {\bf 95},  (1990), 89-114.\\
$\text{[L2]}$ L.Lusztig. {\it Modular representations and quantum groups.}
 Contemp. Math.  \hbox{\bf 82},  (1989), 59-77.\\
$\text{[V]}$
A.Voronov. {\it Semi-infinite homological algebra.} Invent. Math. {\bf 113},
 (1993),  103--146.        \\
$\text{[Xi]}$ Xi Nanhua.  {\it Representations of finite dimensional Hopf
algebras arising from quantum groups.} Preprint,  (1989).\\

\end{document}